\documentclass[10pt]{amsart}
\usepackage{amsfonts}
\usepackage{amsmath}
\usepackage{amsthm}
\usepackage{amssymb}
\usepackage{latexsym}
\usepackage{multicol}
\usepackage{verbatim}

\advance\textwidth by 1.2in \advance\oddsidemargin by -.6in
\advance\evensidemargin by -.6in
\newtheorem{cor}{Corollary}[section]
\newtheorem{lem}[cor]{Lemma}

\newtheorem{prop}[cor]{Proposition}

\theoremstyle{definition}
\newtheorem{defn}{Definition}[section]
\theoremstyle{definition}
\newtheorem{thm}{Theorem}

\newtheorem*{rem}{Remark}

\newenvironment{pf}{\proof}{\endproof}
\newcounter{cnt}
\newenvironment{enumerit}{\begin{list}{{\hfill\rm(\roman{cnt})\hfill}}{%
\settowidth{\labelwidth}{{\rm(iv)}}\leftmargin=\labelwidth%
\advance\leftmargin by
\labelsep\rightmargin=0pt\usecounter{cnt}}}{\end{list}}

\theoremstyle{remark}

\numberwithin{equation}{section} 

\def\dnl{=\hspace{-.2cm}>\hspace{-.2cm}=}
\def\dnr{=\hspace{-.2cm}<\hspace{-.2cm}=}
\def\sn{-\hspace{-.17cm}-}

\begin{document}

\newcommand{\thmref}[1]{Theorem~\ref{#1}}
\newcommand{\secref}[1]{Section~\ref{#1}}
\newcommand{\lemref}[1]{Lemma~\ref{#1}}
\newcommand{\propref}[1]{Proposition~\ref{#1}}
\newcommand{\corref}[1]{Corollary~\ref{#1}}
\newcommand{\remref}[1]{Remark~\ref{#1}}
\newcommand{\defref}[1]{Definition~\ref{#1}}
\newcommand{\er}[1]{(\ref{#1})}
\newcommand{\id}{\operatorname{id}}
\newcommand{\tensor}{\otimes}
\newcommand{\nc}{\newcommand}
\newcommand{\rnc}{\renewcommand}
\newcommand{\qbinom}[2]{\genfrac[]{0pt}0{#1}{#2}}
\nc{\cal}{\mathcal} \nc{\goth}{\mathfrak} \rnc{\bold}{\mathbf}
\renewcommand{\frak}{\mathfrak}
\newcommand{\desc}{\operatorname{desc}}
\newcommand{\Maj}{\operatorname{Maj}}
\renewcommand{\Bbb}{\mathbb}
\nc\bpi{{\mbox{\boldmath $\pi$}}}
\newcommand{\lie}[1]{\mathfrak{#1}}
\makeatletter
\def\section{\def\@secnumfont{\mdseries}\@startsection{section}{1}%
  \z@{.7\linespacing\@plus\linespacing}{.5\linespacing}%
  {\normalfont\scshape\centering}}
\def\subsection{\def\@secnumfont{\bfseries}\@startsection{subsection}{2}%
  {\parindent}{.5\linespacing\@plus.7\linespacing}{-.5em}%
  {\normalfont\bfseries}}
\makeatother
\def\subl#1{\subsection{}\label{#1}}

\nc{\Cal}{\cal} \nc{\Xp}[1]{X^+(#1)} \nc{\Xm}[1]{X^-(#1)}
\nc{\on}{\operatorname} \nc{\ch}{\mbox{ch}} \nc{\Z}{{\bold Z}}
\nc{\J}{{\cal J}} \nc{\C}{{\bold C}} \nc{\Q}{{\bold Q}}
\renewcommand{\P}{{\cal P}}
\nc{\N}{{\Bbb N}} \nc\boa{\bold a} \nc\bob{\bold b} \nc\boc{\bold
c} \nc\bod{\bold d} \nc\boe{\bold e} \nc\bof{\bold f}
\nc\bog{\bold g} \nc\boh{\bold h} \nc\boi{\bold i} \nc\boj{\bold
j} \nc\bok{\bold k} \nc\bol{\bold l} \nc\bom{\bold m}
\nc\bon{\bold n} \nc\boo{\bold o} \nc\bop{\bold p} \nc\boq{\bold
q} \nc\bor{\bold r} \nc\bos{\bold s} \nc\bou{\bold u}
\nc\bov{\bold v} \nc\bow{\bold w} \nc\boz{\bold z}

\nc\ba{\bold A} \nc\bb{\bold B} \nc\bc{\bold C} \nc\bd{\bold D}
\nc\be{\bold E} \nc\bg{\bold G} \nc\bh{\bold H} \nc\bi{\bold I}
\nc\bj{\bold J} \nc\bk{\bold K} \nc\bl{\bold L} \nc\bm{\bold M}
\nc\bn{\bold N} \nc\bo{\bold O} \nc\bp{\bold P} \nc\bq{\bold Q}
\nc\br{\bold R} \nc\bs{\bold S} \nc\bt{\bold T} \nc\bu{\bold U}
\nc\bv{\bold V} \nc\bw{\bold W} \nc\bz{\bold Z} \nc\bx{\bold X}

\title{Spectral characters of finite--dimensional representations of affine 
algebras}

\author{Vyjayanthi Chari and Adriano A. Moura }
\address{Department of Mathematics, University of
California, Riverside, CA 92521.} \email{chari@math.ucr.edu,
adrianoam@math.ucr.edu}

\maketitle

\section*{Introduction}

In this paper we study the category $\cal C$ of
finite--dimensional representations of affine Lie algebras. The
irreducible objects of this category were classified and described
explicitly in  \cite{C},\cite{CPnew}. It was known however that
$\cal C$ was not semisimple. In such a case a natural problem is
to describe the blocks of the category. The blocks of an abelian
category are themselves  abelian subcategories,  each of which
cannot be written as a proper  direct sum of abelian categories
and such that their direct sum is equal to the original category.
Block decompositions of representations of algebras  are often
given by a character, usually a central character, namely a
homomorphism from the center of the algebra to $\bc$, as for
instance in the case of modules from the BGG category $\cal O$ for
a simple Lie algebra. In our case however, the center of the
universal algebra of the affine algebra  acts trivially on all
representations in the category $\cal C$ and the absence of a
suitable notion of character has been an obstacle to determining
the blocks of $\cal C$.

In recent years the study of the corresponding category $\cal C_q$
of modules for quantum affine algebras has been of some interest,
\cite{CPqa}, \cite{CPweyl}, \cite{FM}, \cite{FR}, \cite{Ka2},
\cite{N1}, \cite{N2}. In \cite{EM} the authors defined the notion
of an elliptic character for objects of $\cal C_q$ when $|q|\neq
1$ and showed that for $|q|<1$, the character could be used to
determine the blocks of $\cal C_q$. The original definition of the
elliptic character used
convergence properties of the (non--trivial) action of
the $R$--matrix on the tensor product of finite--dimensional 
representations.
Of course in the $q=1$ case, the action of the $R$--matrix on a
tensor product is trivial. However, the combinatorial part of the
proof given in \cite{EM} suggests that an elliptic character can
be viewed as a function $\chi: E \to \bz^m$ with finite support,
where $E$ is the elliptic curve $\bc^{\times}/q^{2\bz}$ and $m\in
\bn^+$ depends on the underlying simple Lie algebra. This then
motivated our definition when $q=1$ of a spectral character of
$L(\lie g)$ as a function $\bc^\times\to \Gamma$ with finite
support, where $\Gamma$ is the quotient of the weight lattice of
$\lie g$ by the root lattice of $\lie g$.

The other ingredient used in \cite{EM}  to prove that two modules
with the same elliptic character belonged to the same block, was a
result proved in \cite{C},\cite{Ka2} that a suitable
tensor product of irreducible representations was indecomposable
but reducible on certain natural vectors. In the classical case
however, it was known from the work of \cite{CPnew} that a tensor
product of irreducible representations was either irreducible or
completely reducible. However, it was shown in \cite{CPweyl} that
the the tensor product of the irreducible representations of the
quantum affine algebra specialized to indecomposable, but usually
reducible representations of the classical affine algebra.  This
led to the definition of the Weyl modules as   a family of
universal indecomposable modules. The Weyl modules are in general
not well--understood, see \cite{CPweyl},\cite{FL1},\cite{FL2} for
several conjectures about them. However in this paper, we are
still able to identify a large family of quotients of the Weyl
modules, which allow us to effectively use them as a substitute
for the methods of \cite{EM}. Although, we work with the affine
Lie algebra, our results and proofs work for the current algebra,
$\lie g\otimes\bc[t]$, but with  the spectral character being
defined as functions  from  $\bc$ to $\Gamma$ with finite support.

The paper is organized as follows: section 1 is devoted to
preliminaries, section 2 to the definition of the spectral
character and the statement of the main theorem. In Section 3, we
recall the definition of the Weyl modules, give an explicit
realization of certain indecomposable but reducible quotients of
these modules and  the parametrization of the irreducible objects
of $\cal C$. The theorem is proved in the remaining two sections.
We prove that to every indecomposable object of $\cal C$, one can
associate a  spectral character. To do this we  show that if two modules 
$V_j$, $j=1,2$ have distinct spectral characters
then the corresponding $\rm{Ext}^1(V_1,V_2)=0$.
Finally, we prove that any two
modules with the same spectral character must be in the same block
of $\cal C$ and hence we get a parametrization of the blocks of
$\cal C$ analogous to the one in \cite{EM}.

\vskip 12pt

\noindent{\bf Acknowledgments:}

We  thank P.
Etingof for several discussions and for the uniform proof of Proposition 
1.2.
However, in the appendix we give a more elementary but explicit proof of 
this proposition, which is useful
for computations.
The second author
also thanks P. Etingof for
his continued support and
encouragement.

\setcounter{section}{0}

\section{Preliminaries}

Throughout this paper~$\bn$ (respectively, $\bn^+$) denotes the
set of non-negative (respectively, positive) integers.

Let $\frak g$ be a complex finite-dimensional simple
Lie algebra  of rank $n$ with a Cartan subalgebra~$\frak h$. Set
$I=\{1,2,\cdots ,n\}$ and let $\{\alpha_i\,:\,i\in I\}\subset\frak
h^*$ (resp. $\{\omega_i\,:\,i\in I\}\subset\frak h^*$) be the set
of simple roots (resp.  fundamental weights) of $\frak g$ with
respect to $\frak h$. Define a non--degenerate bilinear form $<\ ,
\  >$ on $\lie h$ by $<\omega_i,\alpha_j>=\delta_{ij}$ and let
$h_{i}\in\lie h$ be defined by requiring
$\omega_i(h_{j})=\delta_{ij}$, $i,j\in I$.
We shall assume that the nodes of the Dynkin diagram are numbered
as follows and we let $I_\bullet\subset I$ be the indices of the
shaded nodes in the diagram.

\vspace{.5cm}

{\centerline {\bf Table 1}}
\begin{multicols}{2}
\begin{itemize}

\item $A_n$ : {\large \vspace{-.55cm} 
$$\stackrel{1}{\bullet}\hspace{-.18cm}\sn\hspace{-.18cm}\stackrel{2}{\circ} 
\dots 
\stackrel{\text{n-1}}{\circ}\hspace{-.29cm}\sn\hspace{-.18cm}\stackrel{\text{n}}{\circ}$$}

\item $B_n$ : {\large \vspace{-.55cm} 
$$\stackrel{1}{\circ}\hspace{-.18cm}\sn\hspace{-.18cm}\stackrel{2}{\circ} 
\dots 
\stackrel{\text{n-1}}{\circ}\hspace{-.18cm}\dnl\hspace{-.07cm}\stackrel{\text{n}}{\bullet}$$}

\item $C_n$ : {\large \vspace{-.55cm} 
$$\stackrel{1}{\bullet}\hspace{-.18cm}\sn\hspace{-.18cm}\stackrel{2}{\circ} 
\dots 
\stackrel{\text{n-1}}{\circ}\hspace{-.18cm}\dnr\hspace{-.07cm}\stackrel{\text{n}}{\circ}$$}

\item $D_n$, $n$ odd : {\large \vspace{-.75cm} 
$$\hspace{.95cm}\stackrel{1}{\circ}\hspace{-.18cm}\sn\hspace{-.18cm}\stackrel{2}{\circ} 
\dots 
\stackrel{\text{n-2}}{\circ}\hspace{-.29cm}\sn\hspace{-.28cm}\stackrel{\text{n-1}}{\circ}$$
\vspace{-.82cm}$$\hspace{1.8cm}|$$
\vspace{-.73cm}$$\hspace{2.11cm}\bullet\text{ \small n}$$}

\item $D_n$, $n$ even : {\large \vspace{-.75cm} 
$$\hspace{.95cm}\stackrel{1}{\circ}\hspace{-.18cm}\sn\hspace{-.18cm}\stackrel{2}{\circ} 
\dots 
\stackrel{\text{n-2}}{\circ}\hspace{-.29cm}\sn\hspace{-.28cm}\stackrel{\text{n-1}}{\bullet}$$
\vspace{-.82cm}$$\hspace{1.8cm}|$$
\vspace{-.73cm}$$\hspace{2.11cm}\bullet\text{ \small n}$$}
\\
\item $E_6$  : {\large \vspace{-.55cm}
$$\stackrel{1}{\bullet}\hspace{-.18cm}\sn\hspace{-.18cm}\stackrel{2}{\circ} 
\hspace{-.18cm}\sn\hspace{-.18cm} 
\stackrel{\text{3}}{\circ}\hspace{-.18cm}\sn\hspace{-.18cm}\stackrel{4}{\circ}\hspace{-.18cm}\sn 
\hspace{-.18cm}\stackrel{5}{\circ}$$
\vspace{-.75cm}$$\hspace{-.00cm}|$$
\vspace{-.75cm}$$\hspace{.28cm}\circ\text{ \footnotesize$6$}$$}

\item $E_7$  : {\large \vspace{-.65cm} 
$$\hspace{.5cm}\stackrel{1}{\bullet}\hspace{-.18cm}\sn\hspace{-.18cm}\stackrel{2}{\circ} 
\hspace{-.18cm}\sn\hspace{-.18cm} 
\stackrel{\text{3}}{\circ}\hspace{-.18cm}\sn\hspace{-.18cm} 
\stackrel{4}{\circ}\hspace{-.18cm}\sn 
\hspace{-.18cm}\stackrel{5}{\circ}\hspace{-.18cm} 
\sn\hspace{-.18cm}\stackrel{6}{\circ}$$
\vspace{-.82cm}$$\hspace{.98cm}|$$
\vspace{-.75cm}$$\hspace{1.26cm}\circ\text{ \footnotesize$7$}$$}

\item $E_8$  : {\large \vspace{-.75cm} 
$$\hspace{1cm}\stackrel{1}{\bullet}\hspace{-.18cm}\sn\hspace{-.18cm}\stackrel{2}{\circ} 
\hspace{-.18cm}\sn\hspace{-.18cm} \stackrel{\text{3}}{\circ}\hspace{-.18cm} 
\sn\hspace{-.18cm}\stackrel{4}{\circ}\hspace{-.18cm}\sn 
\hspace{-.18cm}\stackrel{5}{\circ}\hspace{-.18cm} 
\sn\hspace{-.18cm}\stackrel{6}{\circ} 
\hspace{-.18cm}\sn\hspace{-.18cm}\stackrel{7}{\circ}$$
\vspace{-.82cm}$$\hspace{1.95cm}|$$
\vspace{-.75cm}$$\hspace{2.23cm}\circ\text{ \footnotesize$8$}$$}

\item $F_4$ : {\large \vspace{-.55cm} 
$$\hspace{-.4cm}\stackrel{1}{\bullet}\hspace{-.18cm}\sn\hspace{-.18cm}\stackrel{2}{\circ} 
\hspace{-.08cm}\dnr\hspace{-.08cm}\stackrel{3}{\circ} 
\hspace{-.18cm}\sn\hspace{-.18cm}\stackrel{4}{\circ}$$}

\item $G_2$ : {\large \vspace{-.55cm} 
$$\hspace{-1.4cm}\stackrel{1}{\bullet}\hspace{-.12cm}\equiv\hspace{-.18cm}<\hspace{-.18cm}\equiv 
\hspace{-.1cm}\stackrel{2}{\circ}$$}

\end{itemize}

\end{multicols}

Let $R^+$ be the corresponding
set of positive roots and denote by~$\theta$ the highest root
of~$\lie g$. As usual, $Q$ (respectively, $P$) denotes the root
(respectively, weight) lattice of $\frak g$ and we let
$\Gamma=P/Q$. It is known that, \cite{H},
\begin{eqnarray*} &\Gamma\cong& \bz_{n+1},\ \lie g \ {\text{of type}}\ \
A_n,\ \ \   \Gamma\cong  \bz_2,\ \ \lie g \  {\text{of type}}\ \
B_n, C_n, E_7\\ &\Gamma\cong & \bz_4,\ \  \lie g\ {\text{of
type}}\ \ D_{2m+1},\ \ \ \  \Gamma\cong \bz_2\times \bz_2,\ \ \lie
g \ {\text{of type}}\ \ D_{2m},\\ &\Gamma\cong &\bz_3,\ \  \lie g\
{\text{of type}}\  E_6,\ \ \ \ \ \  \Gamma\cong 0,\ \ \lie g\
{\text{of type}}\ \ E_8,F_4,G_2
\end{eqnarray*}
The group $\Gamma$ is generated by the images of the elements
$\{\omega_i:i\in I_\bullet\}$ and hence any $\gamma\in \Gamma$
defines a unique element $\lambda_\gamma=\sum_{i=1}^n
r_i\omega_i\in P^+$ where $r_i\in\bz$ are the minimal
non--negative integers such that $\lambda_\gamma$ is a
representative of $\gamma$. In particular, $r_i=0$ if $i\notin
I_\bullet$.

Let $W$ be the Weyl group of $\lie g$ and assume that $w_0$ is the
longest element of $W$. The group $W$ acts on $\lie h^*$ and
preserves the root and weight lattice.
Let~$P^+=\sum_{i\in I}\bn\omega_i$ be the set of
dominant integral  weights and set~$Q^+=\sum_{i\in I}\bn\alpha_i$.
We shall assume  that $ P$ has the usual partial  ordering, given
$\lambda,\mu\in P$ we say that $\lambda\ge \mu$ if $\lambda-\mu\in
Q^+$. For $\alpha\in R^+$, let $\lie g_{\pm\alpha}$ denote the
corresponding root spaces, and fix elements $x_\alpha^\pm\in\frak
g_{\pm\alpha}$, $h_\alpha\in\lie h$,  such that they span a
subalgebra of $\lie g$ which is isomorphic to $sl_2$. For $i\in
I$, set $x_i^\pm=x_{\alpha_i}^\pm$, $h_{\alpha_i}=h_i$.  Set
$\frak n^\pm=\oplus_{\alpha\in R^+}\lie g_{\pm\alpha}.$ Given any
Lie algebra $\lie a$, let $L(\lie a)=\lie a\otimes \bc[t,t^{-1}]$
be the loop algebra associated with~$\frak a$ and let $\bu(\lie
a)$ be the universal enveloping algebra of $\lie a$. Clearly, we
have
\begin{equation*}\lie g=\lie n^+ \oplus\lie h\oplus\lie n^-, \ \
L(\lie g)=L(\lie n^+) \oplus L(\lie h) \oplus L( \lie n^-),
\end{equation*} and a corresponding decomposition
\begin{equation*} \bu(\lie g)=\bu(\lie n^-)\bu(\lie h)\bu(\lie
n^+),\ \  \bu(L(\lie g))=\bu(L(\lie n^-))\bu(L(\lie h))\bu(L(\lie
n^+)).\end{equation*} Given $a\in\bc^\times$, let $ev_a:L(\lie
g)\to \lie g$ be the evaluation homomorphism, $ev_a(x\otimes
t^n)=a^nx$.

  For $\lambda\in P^+$, let $V(\lambda)$ be
the irreducible finite--dimensional $\lie g$--module with highest
weight $\lambda$ and highest weight vector $v_\lambda$. Thus,
$V(\lambda)$ is generated by $v_\lambda$ as   a $\lie g$--module
with defining relations:
\begin{equation*}\lie n^+v_\lambda=0,\ \
hv_\lambda=\lambda(h)v_\lambda,\ \
(x_i^-)^{\lambda(h_i)+1}v_\lambda =0,\ \ \forall\ h\in\lie h, i\in
I .\end{equation*}
Let $V$ be a finite--dimensional representation of $\lie g$. Then we can 
write
$V$ as a direct sum \begin{equation}\label{wtspace}V=\bigoplus_{\mu\in P} 
V_\mu,\ \ \ \  V_\mu=\{v\in
v: h v=\mu(h)v,\ \ \forall\ h\in\lie h\},\end{equation} and set $$\text{wt}
(V)=\{\mu\in P: V_\mu\ne 0\}.$$
The following result is well--known, \cite{H}.
\begin{prop}\label{winv} Let $V$ be a finite---dimensional representation of
$\lie g$.
\begin{enumerit}
\item[(i)] For all  $w\in W$, $\mu\in P$, we have
$\text{dim}(V_\mu)=\text{dim}(V_{w\mu})$.
\item[(ii)]
The module  $V$ is isomorphic to a direct sum of $\lie g$--modules
of type $V(\lambda)$, $\lambda\in P^+$.
\item[(iii)] Let $V(\lambda)^*$ be the representation of $\lie g$ which is  
dual
to $V(\lambda)$. Then $$V(\lambda)^*\cong V(-w_0\lambda).$$
\hfill\qedsymbol
\end{enumerit}\end{prop}

The following proposition is crucial for the proof of the main
theorem.

\begin{prop}\label{shaded}
Let $\mu,\lambda\in P^+$ be such that $\lambda-\mu\in Q$. Then,
there exists a sequence of weights $\mu_l \in P^+$, $l=0,\cdots,
m$, with

\begin{enumerit}\item[(i)]  $\mu_0 = \mu$,   $\mu_m=\lambda$,
  and
\item[(ii)] ${\rm Hom}_{\lie g}(\lie g\otimes V(\mu_l), V(\mu_{l+1}))\ne 0,\
\ \forall\ 0\le l\le m$.
\end{enumerit}
\end{prop}

\begin{pf} Consider the module $V(\lambda)\otimes V(\mu)^*$. Since 
$\lambda-\mu\in Q$
it follows that $\lambda-w_0\mu\in Q$. In particular this means that if 
$V(\nu)$ is an irreducible summand
of $V(\lambda)\otimes V(\mu)^*$ then $\nu\in Q^+\cap P^+$. It follows that 
$V(\nu)_0\ne 0$.
This implies by a result of Kostant \cite{K,D}, 
that there exists $m\ge 0$ such that $\text{Hom}_{\lie g}(S^m(\lie g), 
V(\nu))\ne 0$. It follows that
$\text{Hom}_{\lie g}(S^m(\lie g)\otimes V(\mu),  V(\lambda))\ne 0$

We now  proceed by induction on $m$.
If $m=1$, we are done for then $\mu_0=\mu$. Otherwise there must  exist
$\mu_{m-1}\in P^+$ with  $\text{Hom}_{\lie  g}(V(\mu_{m-1}), \lie 
g^{\otimes(m-1)}\otimes V)\ne 0 $
such that  $\text{Hom}_{\lie g}(\lie g\otimes V(\mu_{m-1}), U)\neq 0$.  
Since the
category of finite--dimensional representations of $\lie g$ is semisimple, 
we see also that
$\text{Hom}_{\lie g}(\lie g^{\otimes (m-1)}\otimes V, V(\mu_{m-1}))\neq 0$. 
But now we are done by the inductive hypothesis.
\end{pf}

\begin{rem}
In the appendix   we construct the
sequence $\mu_1,\dots, \mu_m$ in the special case when $\mu=\lambda_\gamma$
with the properties stated above. In particular, this gives a different and 
perhaps more elementary proof of
this proposition.

\end{rem}

\section{ Spectral characters and the block decomposition of $\cal C$.}

Let $\Xi$ be the set of all functions $\chi: \bc^\times
\to\Gamma$ with finite support. Clearly addition of functions
defines a group structure on $\Xi$. Given $\lambda\in P^+$,
$a\in\bc^\times$, let $\chi_{\lambda, a}\in \Xi$ be defined by
\begin{equation*}\chi_{\lambda,a}(z)=\delta_a(z)\overline{\lambda},\end{equation*}
where $\overline\lambda$ is the image of $\lambda$ in $\Gamma$ and
$\delta_a(z)$ is the characteristic function of $a\in\bc^\times$.
We denote by $\cal{P}$ the space of $n$--tuples of polynomials
with constant term one. Coordinatewise multiplication defines the
structure of monoid on $\cal P$. Given $\lambda\in P^+$,
$a\in\bc^\times$, define $\bpi_{\lambda,a}=(\pi_1,\cdots
,\pi_n)\in\cal P$, by
\begin{equation*}\pi_i=(1-au)^{\lambda(h_i)},\ \ \
1\le i\le n\end{equation*}
Any $\bpi=(\pi_1,\cdots ,\pi_n)\in\cal P$
can be written uniquely as a product
\begin{equation}\label{factorpi} 
\bpi=\prod_{j=1}^r\bpi_{\lambda_j,a_j},\end{equation}  where\hfill
\begin{enumerit}
\item[(i)]$\{a_j^{-1}
:1\le j\le r\}$ is the set of distinct roots of
$\prod_{i=1}^n\pi_i$,
\item[(ii)] $\lambda_j=\sum_{k=1}^n m_{kj}\omega_k\in P^+$, and $m_{kj}$ is 
the multiplicity
with which $a_j^{-1}$ occurs as a root of $\pi_k$.
\end{enumerit} Define $\bpi^*\in\cal P$ by \begin{equation*}
\bpi^*=\prod_{j=1}^r\bpi_{-w_0\lambda_j,a_j},\end{equation*} where
we recall that $w_0$ is the longest element of the Weyl group of
$\lie g$.
Given $\bpi\in\cal P$ define $\chi_\bpi\in\Xi$  by
$$\chi_\bpi=\sum_{j=1}^r\chi_{\lambda_j,a_j},$$ where
$\lambda_i,a_i$ are as in \eqref{factorpi}. Obviously,
$$\chi_{\bpi\bpi'}=\chi_\bpi+\chi_{\bpi'},$$ for all
$\bpi,\bpi'\in\cal P$.

To state our main result we need to recall the
parametrization of irreducible finite--dimensional modules of
affine Lie algebras, \cite{C},\cite{CPweyl} and also the defintion
of blocks in an abelian category.
\begin{prop} There exists a bijective correspondence between the
isomorphism classes of  irreducible finite--dimensional
representations of affine Lie algebras and elements of $\cal P$.
\end{prop} We denote by $V(\bpi)$ an element of the isomorphism
class corresponding to $\bpi$.

\begin{defn} We say that a module $V\in\cal C$ has spectral
character $\chi\in\Xi$ if $\chi=\chi_\bpi$ for every irreducible
component $V(\bpi)$ of $V$.  Let   $\cal C_\chi$ be the abelian
subcategory consisting of all modules $V\in\cal C$ with spectral
character $\chi$.
\end{defn}

Let  $\cal C=\cal
C_{fin}(L(\lie g))$ be the category of finite--dimensional
representations of $L(\lie g)$. This category is not semisimple,
i.e., there exist indecomposable reducible $L(\lie g)$--modules in
$\cal C$. However  $\cal C$ is an abelian tensor category and
every  object in $\cal C$ has a Jordan--Holder series of finite
length. This means that $\cal C$ has a block decomposition which
is obtained as follows.
\begin{defn} Say that
two indecomposable  objects $U,V\in\cal C$ are linked if there do
not exist abelian subcategories $\cal C_k$, $k=1,2$ such that
$\cal C=\cal C_1\oplus \cal C_2$ with $U\in\cal C_1$, $V\in\cal
C_2$. If $U$ and $V$ are decomposable then we say that they are
linked if every indecomposable summand of $U$ is linked to every
indecomposable summand of $V$. \end{defn}
  This  defines an
equivalence relation on $\cal C$ and  a  block of $\cal C$  is an
equivalence class for this relation, clearly $\cal C$ is a direct
sum of blocks. The following lemma is trivially established.
\begin{lem}
Two indecomposable modules  $V_1$ and $V_2$ are linked
iff they contain  submodules $U_k\subset V_k$, $k=1,2$ such that
$U_1$ is linked to $U_2$.
\end{lem}

The main result of the paper is the following.
\begin{thm}\label{main} We have $$\cal C=\bigoplus\limits_{\chi\in\Xi}\cal 
C_\chi.$$
Moreover each $\cal C_\chi$ is a block. Equivalently, the blocks
of $\cal C$ are in bijective correspondence with $\Xi$.
\end{thm}

The theorem is obviously a consequence of the next two
propositions.
\begin{prop}\label{irrlink} Any two irreducible modules in $\cal
C_\chi$, $\chi\in\Xi$, are linked.
\end{prop}
\begin{prop}\label{indecomp} Every  indecomposable $L(\lie g)$--module has a
spectral character.
\end{prop}

We prove these propositions in Sections 4 and 5 of the paper
respectively. We shall need several results on a certain family of
indecomposable but generally reducible modules for $L(\lie g)$,
the so--called Weyl modules, this is done  in the next section.

We conclude this section with an equivalent definition of   linked modules 
and
with some general results on Jordan--Holder series.
\begin{defn} Let $U,V\in\cal C$ be indecomposable  $L(\lie g)$--modules.
We say that $U$ is strongly linked to $V$ is there exists  $L(\lie
g)$--modules $U_1,\cdots ,U_\ell$, with $U_1=U$, $U_\ell=V$ and
either ${\rm{Hom}}_{L(\lie g)}(U_k,U_{k+1})\neq 0$ or
${\rm{Hom}_{L(\lie g)}}(U_{k+1},U_{k})\neq 0$ for all $1\le k\le
\ell $. We extend this to all of $\cal C$ by saying that two
modules $U$ and $V$ are strongly linked iff every indecomposable
component of $U$ is strongly linked to every indecomposable
component of $V$.
\end{defn}

It is clear that the notion of strongly linked defines an
equivalence relation on $\cal C$ which induces a decomposition of
$\cal C$ into a direct sum of abelian categories.  If two  modules
$U$ and $V$ are strongly linked then they must be linked. For
otherwise, suppose that $U$ and $V$ belong to different blocks. It
suffices to consider the case $\text{Hom}_{L(\lie g)}(U,V)\ne 0$.
This means that $U$ and $V$ have an irreducible constituent say
$M$ in common. Then, since each block is an abelian subcategory,
$M$ must belong to both blocks which is a contradiction.
Conversely, suppose that $U$ and $V$ are linked but not strongly
linked. Then, there is obviously a splitting of $\cal C$ into
abelian subcategories coming from the strong linking, such that
$U$ and $V$ belong to different subcategories. We have proved the
following:
\begin{lem} Two modules $U$ and $V$ are linked iff they are
strongly linked.\end{lem}

\begin{lem} \label{stronglink} Suppose that $U\in\cal C_{\chi_1}$ and 
$V\in\cal C_{\chi_2}$ are  strongly
linked. Then $\chi_1=\chi_2$. \end{lem}
\begin{pf} It suffices to check this when $\text{Hom}_{L(\lie
g)}(U,V)\ne 0$. But this means that $U$ and $V$ have an
irreducible constituent say $M$ in common and hence
$\chi_1=\chi_2$.
\end{pf}

We shall make use of the following simple proposition repeatedly
without further mention.
\begin{prop}\hfill
\begin{enumerit}

\item[(i)] Any sequence $0\subset V_1\cdots \subset V_k\subset V$
of $L(\lie g)$--modules in $\cal C$ can be refined to a
Jordan--Holder series of $V$.

\item[(ii)] Suppose that $0\subset U_1\subset\cdots\subset U_r=U$
and $0\subset V_1\cdots \subset V_s\subset V$ are Jordan--Holder
series for modules $U$, $V$ in $\cal C$. Then the irreducible
constituents of $U\otimes V$ occur as constituents of $U_k\otimes
V_\ell$ for some $1\le k\le r$ and $1\le \ell\le s$.

\item[(iii)]\label{tstronglink}
Suppose that $U_k\in\cal C$, $1\le k\le 3$ and that $U_1$ and
$U_2$ are strongly linked.  Then $U_1\otimes U_3$ is strongly
linked to $U_2\otimes U_3$.

\end{enumerit}
\end{prop}

\section{Weyl modules} In this section we recall from \cite{CPweyl}
the definition and some results
on Weyl modules. We also study further properties of these
modules.

Let $V\in\cal C$. Regarding $V$ as a finite--dimensional module for $\lie 
g$,  we can
write $V$ as a direct sum as in Section 1, \eqref{wtspace}
$V=\oplus_{\mu\in P} V_\mu$. Let $\text{wt}\ (V)$ be the set of
weights of $V$. Notice that $L(\lie h)V_\mu\subset V_\mu$. Since
$L(\lie h)$ is an abelian Lie algebra, we get a further
decomposition $$V_\mu=\bigoplus_{\bod\in L(\lie h)^*}
V_\mu^\bod,$$ where $$V_\mu^\bod=\{v\in V_\mu: \left(h\otimes
t^k-\bod(h\otimes t^k)\right)^r v=0, \ \ \forall \ r \ge
r(h,k)>>0\},$$ are the generalized eigenspaces for the action of
$L(\lie h)$ on $V_\mu$.  Clearly if $U,V\in\cal C$, then any
$L(\lie g)$ homomorphism from $U$ to $V$ maps $U_\mu^\bod$ to
$V_\mu^\bod$. Since $V_\mu$ is finite--dimensional we see that if
$V_\mu^\bod\ne 0$ then there exists $0\ne v\in V_\mu^\bod$ such
that $$(h\otimes t^k) v=\bod(h\otimes t^k)v,\ \ h\in\lie h,\ \
k\in\bz.$$  We say that $\bod$ is of type $\bpi\in\cal P$,  if
$$\bod(h\otimes t^k)=(\sum_{j=1}^r\lambda_j(h)a_j^k),$$ where
$\lambda_j\in P^+$ and $a_j\in\bc^\times$ are as in
\eqref{factorpi} and we denote the corresponding generalized
eigenspace by $V_\mu^\bpi$.

\begin{defn} Given an $n$--tuple of polynomials with constant term 1,
we denote by $W(\bpi)$ the $L(\lie
g)$--module generated by an element $w_\bpi$ and the following
relations:
\begin{equation}\label{relweyl} L(\lie n^+)w_\bpi=0,\ \ \ (h\otimes
t^k) w_\bpi=(\sum_{j=1}^r\lambda_j(h)a_j^k)w_\bpi,\ \ \
(x_{i}^-\otimes t^ \ell)^{\sum_{j=1}^r\lambda_j(h_i)+1}w_\bpi =0,
\end{equation}
for all $i\in I$,  $k,\ell\in\bz$, $\alpha\in R^+$, $h\in\lie h$
and where we assume that $\bpi$ is written as in \eqref{factorpi}.
Set
$\lambda_\bpi=\sum_{j=1}^r\lambda_j$.
\end{defn}
The following properties of $W(\bpi)$ are standard and  easily
established:
\begin{lem}\label{stweyl} With the notation as above, we have,\hfill
\begin{enumerit}
\item[(i)] $W(\bpi)=\bu(L(\lie n^-)) w_\bpi$ and so
$\text{wt}(W(\bpi))\subset \lambda_\bpi-Q^+$.
\item[(ii)]
$\text{dim}\  W(\bpi)_{\lambda_\bpi}=1$, and so
$W(\bpi)^\bpi_{\lambda_\bpi}=W_{\lambda_\bpi}$.
\item[(iii)] Let $V$ be
any finite--dimensional $L(\lie g)$--module generated by an
element $v\in V$ satisfying
$$ L(\lie n^+)v=0,\ \ \ (h\otimes
t^k) v=(\sum_{j=1}^r\lambda_j(h)a_j^k)v. $$
Then $V$ is a quotient of $W(\bpi)$.
\item[(iv)]$W(\bpi)$
is an indecomposable $L(\lie g)$--module with a unique irreducible
quotient.
\end{enumerit}\hfill\qedsymbol
\end{lem}

Let $\beta_1,\cdots ,\beta_N$ be an enumeration of the
elements of $R^+$. Given $r\in\bz$, set
$x^-_{\beta_j,r}=x_{\beta_j}\otimes t^r$. The following result was
proved in \cite{CPweyl}.
\begin{thm}\label{weylprop} Let $\bpi$ be an $n$--tuple of polynomials with
constant term one and assume that  $\bpi$ has a factorization  as
in \eqref{factorpi}.
\begin{enumerit}
\item[(i)] The $L(\lie g)$--module $W(\bpi)$ is spanned by
monomials of the form
$$x^-_{\beta_{j_1},r_1}x^-_{\beta_{j_2},r_2}\cdots
x^-_{\beta_{j_\ell},r_\ell}w\bpi$$ where $\ell\in\bn^+$, $j_1\le
j_2\le\cdots\le  j_\ell$ and  $0\le r_k<
\lambda_\bpi(h_{\beta_{j_k}})$ for all $1\le k\le \ell$. In
particular,  $\text{dim}\ W(\bpi)<\infty$.
\item[(ii)] As
$L(\lie g)$--modules,
$$W(\bpi)\cong W(\bpi_{\lambda_1,a_1})\otimes\cdots\otimes
W(\bpi_{\lambda_r, a_r}).$$
\end{enumerit}\hfill\qedsymbol
\end{thm}
We can now elaborate on the parametrization of the irreducible
finite--dimensional modules stated in Section 2 of this paper.
\begin{prop}\label{irrfactor}  The irreducible finite--dimensional $L(\lie
g)$--module $V(\bpi)$ is the irreducible quotient of $W(\bpi)$ and
we have $$V(\bpi)\cong V(\bpi_{\lambda_1,a_1})\otimes\cdots\otimes
V(\bpi_{\lambda_r, a_r}).$$  Further, the module $V(\bpi_{\lambda,
a})$ is the $L(\lie g)$--module obtained by pulling back the $\lie
g$--module $ V(\lambda)$, by the evaluation homomorphism
$\text{ev}_a:L(\lie g)\to \lie g$. Finally  as  $L(\lie
g)$--modules we have
$$V(\bpi)^*\cong V(\bpi^*).$$\hfill\qedsymbol
\end{prop}

The structure of $W(\bpi)$ is not well--understood in general,
although it is known that $W(\bpi)$ is in general not isomorphic
to $V(\bpi)$, a necessary and sufficient condition for $W(\bpi)$ to be
isomorphic to $V(\bpi)$ can be found in \cite{CPweyl}.  In what
follows, we establishe further properties of the Weyl modules
which we need in this paper, and also identify natural
indecomposable reducible  quotients  of $W(\bpi)$.

\begin{prop}\label{specprop} Let $\lambda=\sum_{i=1}^nr_i\omega_i\in P^+$,
$a\in\bc^\times$.
\begin{enumerit}
\item[(i)] For all $\alpha\in R^+$ we have $$(x_\alpha^-\otimes
(t-a)^{\lambda(h_\alpha)})w_{\bpi_{\lambda,a}} =0.$$ In particular
$W(\bpi_{\lambda, a})$ is spanned by elements of the
form$$(x^-_{\beta_{j_1}}\otimes(t-a)^{r_1})(x^-_{\beta_{j_2}}\otimes
(t-a)^{r_2})\cdots (x^-_{\beta_{j_\ell}}\otimes
(t-a)^{r_\ell})w\bpi_{\lambda,a} $$ where $\ell\in\bn^+$, $j_1\le
j_2\le\cdots\le j_\ell$ and  $0\le r_k<
\lambda_\bpi(h_{\beta_{j_k}})$ for all $1\le k\le \ell$.
\item[(ii)] For all $h\in \lie h$, $k\in\bz$, $\mu\in P$ and $w\in
W(\bpi_{\lambda,a})_\mu$, we have,
$$(h\otimes (t^k-a^k))^r w=0,\ \ \forall \ r>>0.$$
\item[(iii)] There exists a bijective correspondence
between irreducible $\lie g$--submodules of $W(\bpi_{\lambda, a})$
and the irreducible $L(\lie g)$--constituents of $W(\bpi_{\lambda,
a})$.
\end{enumerit}
\end{prop}

\begin{pf}
The relation 
$(x_\alpha^-\otimes(t-a)^{\lambda(h_\alpha)})w_{\bpi_{\lambda,a}} =0$ was 
proved in section 6 of \cite{CPweyl}. This immediately implies the second 
assertion of $(i)$. To prove $(ii)$ one just uses commutation relations once 
we know from $(i)$ that $w = 
(x^-_{\beta_{j_1}}\otimes(t-a)^{r_1})(x^-_{\beta_{j_2}}\otimes
(t-a)^{r_2})\cdots (x^-_{\beta_{j_\ell}}\otimes
(t-a)^{r_\ell})w\bpi_{\lambda,a} $.

To prove $(iii)$ first notice that, from $(ii)$, it follows that the 
irreducible constituents of
$W(\bpi_{\lambda,a})$ are all of the form $V(\bpi_{\mu, a})$ for some 
$\mu\in P^+$. Then, since $V(\bpi_{\mu, a})\cong_{\lie g} V(\mu)$, it 
follows that all $\lie g$-constituents of $W(\bpi_{\mu, a})$ must also be 
$L(\lie g)$-constituents with the same multiplicity.
\end{pf}

We now prove,
\begin{prop}\label{derv} Let $\lambda, \mu\in P^+$. Assume that  there 
exists a non--zero homomorphism $p:\lie
g\otimes V(\lambda)\to  V(\mu)$ of $\lie g$--modules.  The
following formulas define an action of $L(\lie g)$--module  on
$V(\lambda)\oplus V(\mu)$: $$xt^r(v,w)= (a^rxv,
a^rxw+ra^{r-1}p(x\otimes v)),$$ where $x\in\lie g$, $r\in\bz$,
$v\in V(\lambda)$ and $w\in V(\mu)$. Denoting this module by
$V(\lambda, \mu,a)$, we see that $$0\to V(\bpi_{\mu,a})\to
V(\lambda,\mu,a)\to V(\bpi_{\lambda,a})\to 0$$ is a non--split
short exact sequence of $L(\lie g)$--modules. Finally, if $\lambda>\mu$,  
there
exists a canonical surjective homomorphism of $L(\lie g
)$--modules $W(\bpi_{\lambda, a})\to V(\lambda ,\mu, a)$.
\end{prop}
\begin{pf} To check that the formulas give a $L(\lie g)$--module
structure is a straightforward verification. Since $L(\lie
g)V(\mu)\subset V(\mu)$ it follows that $V(\bpi_{\mu,a})$ is a
$L(\lie g)$--submodule  of $V(\lambda,\mu,a)$. Since $p:\lie
g\otimes V(\lambda)\to  V(\mu)$ is non--zero, it follows that the
module $V(\lambda,\mu,a)$ is indecomposable and we have the
desired short exact sequence of $L(\lie g)$--modules. Note that
if $\lambda>\mu$, we have $$L(\lie n^+)( v_\lambda,\  0)
=0,\ \ h\otimes t^k (v_\lambda,\ 0)=(a^kv_\lambda,\ 0).$$ Also,
since $V(\lambda)=\bu(\lie g)v_\lambda$ we see that $(V(\lambda),
0) \subset \bu(L(\lie g))v_\lambda$,  and hence it follows that
$V(\lambda,\mu,a)=\bu(L(\lie g))v_\lambda$. But now Lemma
\ref{stweyl} (iii) implies that $V(\lambda, \mu,a)$ must be a
quotient of $W(\bpi_{\lambda,a})$.
\end{pf}

\vskip 12pt

\noindent{\bf{Remark}} One can view the modules $V(\lambda,\mu,a)$
as generalizations of the modules $V(\bpi_{\lambda,a})$ as
follows. Thus, while  $V(\bpi_{\lambda, a})$ is a module for
$L(\lie g)$ on which $x\otimes (f-f(a))$ acts trivially for all
$f\in \bc[t,t^{-1}]$ and $x\in\lie g$, the modules $V(\lambda,
\mu, a)$ are modules on which  $x\otimes (f-(t-a)f'(a)-f(a))$ acts
trivially for all $f\in \bc[t,t^{-1}]$, where $f'$ is the first
derivative of $f$ with respect to $t$. This obviously gives rise
to the natural question of describing  modules of $L(\lie g)$ such
that $x\otimes \left(f-\sum_{j=0}^n (t-a)^jf^{(j)}(a)\right)$ acts
trivially for all $x\in\lie g$ and $f\in\bc[t,t^{-1}]$.

\vskip 24pt

\section{Proof of Proposition \ref{irrlink}}
We begin with the following lemma.

\begin{lem}\label{specprop1}\hfill
\begin{enumerit}
\item[(i)]  Assume that $\lambda,\mu\in P^+$ and that there exists a
non--zero homomorphism $p:\lie g\otimes V(\lambda)\to  V(\mu)$ of
$\lie g$--modules. Then the modules $V(\bpi_{\lambda,a})$ and
$V(\bpi_{\mu,a})$ are strongly linked.
\item[(ii)] Let $\gamma\in\Gamma$ be such that
$\lambda=\lambda_\gamma\mod Q$. Then, $V(\bpi_{\lambda,a})$ and
$V(\bpi_{\lambda_\gamma,a})$ are strongly linked.
\end{enumerit}
\end{lem}
\begin{pf} The first part of the Lemma is immediate from
Proposition \ref{derv}. The second part is now  immediate from {\it(i)} and
Proposition \ref{shaded}.
\end{pf}

\begin{prop}\label{tensor}  Let $V(\bpi_k)\in\cal C_{\chi_k}$ for some 
$\chi_k\in\Xi$,
$k=1,2$. Then $V(\bpi_1)\otimes V(\bpi_2)\in\cal
C_{\chi_1+\chi_2}$
\end{prop}
\begin{pf}  By Proposition \ref{irrfactor} we can write  
$V(\bpi_1)=\otimes_{j=1}^k
V(\bpi_{\lambda_j,a_j})$ with $a_j\ne a_l$ for all $1\le l\ne j\le k$ and
$\lambda_1,\cdots ,\lambda_k\in P^+$. Similarly write
$V(\bpi_2)=\otimes_{j=1}^\ell V(\bpi_{\mu_j,b_j})$. We proceed
by induction on the cardinality of $S$, where
$$S=\{a_1,\cdots ,a_k\}\cap\{b_1,\cdots
,b_\ell\}.$$ If $S$ is empty then $V\otimes U$ is irreducible and
the result is clear. Suppose then that $S\ne \emptyset$ and assume
without loss of generality that $a_1=b_1$. Write
$$V(\bpi_{\lambda_1,a_1})\otimes
V(\bpi_{\mu_1,a_1})=\bigoplus\limits_{\nu\in P^+} m_\nu
V(\bpi_{\nu,a_1})$$ where $m_\nu$ is the multiplicity with which
$V(\nu)$ occurs inside the tensor product of the $\lie g$--modules
$V(\lambda_1)\otimes V(\mu_1)$. Since $\lambda+\mu-\nu\in Q^+$, it
follows from the definition of spectral characters that
$\chi_{\bpi_{\nu,a_1}} =
\chi_{\bpi_{\lambda_1,a_1}}+\chi_{\bpi_{\mu_1,a_1}}$. The
inductive step follows by noting that, $$V(\bpi_1)\otimes
V(\bpi_2)=\left(\bigoplus\limits_{\nu} m_\nu
V(\bpi_{\nu,a_1})\right) \bigotimes\limits_{s=2}^k
V(\bpi_{\lambda_s,a_s})\bigotimes\limits_{j=2}^\ell
V(\bpi_{\mu_j,b_j})$$
\end{pf}

\begin{cor}\label{tensor1} \hfill
\begin{enumerit}
\item[(i)] For all $\chi_k\in\Xi$, $k=1,2$, we have $$\cal
C_{\chi_1}\otimes\cal C_{\chi_2}\subset\cal C_{\chi_1+\chi_2}.$$
\item[(ii)]  Let $V\in\cal C_\chi$, then
$V^*\in\cal C_{-\chi}$.
\end{enumerit}
\end{cor}
\begin{pf}  Let $V_k\in\cal C_{\chi_k}$, $k=1,2$, Since every irreducible 
constituent of
$V_k$ is in $\cal C_{\chi_k}$ part {\it (i)} is immediate from the
proposition.  For part {\it (ii)}, Suppose that $V^*\in\cal C_{\chi'}$
for some $\chi'\in\Xi$. Since $V\otimes V^*$ contains the trivial
representation of $L(\lie g)$ it follows that $V\otimes V^*\in\cal
C_0$ and the lemma is proved.
\end{pf}

\noindent{\it Proof of Proposition \ref{irrlink}.} Suppose that
$V(\bpi_\ell)$, $\ell=1,2$  are irreducible $L(\lie g)$ modules
with the same spectral character $\chi$. By Proposition
\ref{tensor},  there exist $\lambda_{1\ell},\cdots
,\lambda_{s,\ell}\in P^+$, $\ell=1,2$ and $a_1,\cdots
,a_s\in\bc^\times$ such that $\lambda_{1j}-\lambda_{2j}\in Q$ and
$$\bpi_\ell=\prod_{j=1}^s\bpi_{\lambda_{j\ell,a_j}}.$$  If $s=1$,
then the proposition follows from Lemma \ref{specprop1}. If
$\chi=\sum_{j=1}^s\chi_{\lambda_j,a_j}$, then it follows from
Proposition \ref{tstronglink} and Lemma \ref{specprop1} that
$V(\bpi_\ell)$ is strongly linked to $\otimes_{j=1}^s
V(\bpi_{\overline{\lambda}_{\ell,j},a_j})$. The result follows.

\hfill\qedsymbol

\vskip 12pt

\section{Proof of Proposition \ref{indecomp}}
We begin with,
\begin{lem}{\label{oneroot}} We have $W(\bpi)\in\cal C_{\chi_\bpi}$.
\end{lem}
\begin{pf} In view of Corollary \ref{tensor1}, it suffices to prove the 
lemma when
$\bpi=\bpi_{\lambda,a}$.  It follows from Proposition
\ref{specprop} that every irreducible component of $W(\bpi)$ is of
the form $\chi_{\mu,a}$ for some $\mu\in\lambda- Q^+$. The result
is now immediate.
\end{pf}

\begin{lem}\label{nullext}
\begin{enumerit}
\item[(i)]
Let $U\in\cal C_\chi$. Let   $\bpi_0 \in \cal P$
be such that $\chi\ne \chi_{\bpi_0}$.
Then $\rm{Ext}_{L(\lie g)}^1(U,V(\bpi_0)) = 0$.
\item[(ii)]
Assume that  $V_j\in\cal C_{\chi_j}$, $j=1,2$ and that $\chi_1\ne \chi_2$.
Then $\rm{Ext}_{L(\lie g)}^1(V_1,V_2) = 0$.
\end{enumerit}
\end{lem}

\begin{pf} Since $\rm{Ext}^1$ preserves direct sums, to prove (i) it 
suffices to consider the case
when $U$ is indecomposable.
Consider an extension,
\begin{equation*}
0\to V(\bpi_0)\to V\to U \to 0
\end{equation*}
We prove  by induction on the length of $U$ that the extension is trivial.
Suppose first that  $U=V(\bpi)$ for some $\bpi\in\cal P$ and that 
$\chi_{\bpi}\ne\chi_{\bpi_0}$.
Then, either
\begin{enumerit}
\item[(i)] $\lambda_{\bpi}<\lambda_{\bpi_0}$,
\item[(ii)] $\lambda_{\bpi_0}-\lambda_{\bpi} \notin (Q^+-\{0\})$.
\end{enumerit}

Since dualizing the exact sequence above takes us from (i) to (ii),
we can assume without loss of generality that we are in  case (ii). This 
implies immediately that
$$L(\lie n^+) V_{\lambda_\bpi}=0, $$
since $\text{wt}(V(\bpi_0))\subset\lambda_{\bpi_0}-Q^+.$  On the other hand 
since
$V_{\lambda_\bpi}$ maps onto $V(\bpi)_{\lambda_\bpi}$ we see that
$\text{dim}V_{\lambda_{\bpi}}^{\bpi}\ne 0.$ Thus there exists
an element  $0\ne v\in V_{\lambda_\bpi}$
which is a common eigenvector for the action of  $L(\lie h)$ with eigenvalue 
$\bpi$. Since $V$ has length two it follows that either
 $V= \bu(L(\lie g))v$ or that  $\bu(L(\lie g))v=V(\bpi_0)$.
But, the submodule $\bu(L(\lie g))v$ of $V$ is a quotient of 
$W(\bpi)$ and hence has spectral
character $\chi_\bpi$. Since $\chi_\bpi\ne \chi_{\bpi_0}$, we get 
$V(\bpi_0)\cap \bu(L(\lie g))v=0$.
Hence 
$$V\cong V(\bpi_0)\oplus \bu(L(\lie g))v.$$
This proves that induction begins.

Now assume that $U$ is indecomposable but reducible and that we know the 
result for all modules of length
strictly less than that of $U$. Let $U_1$ be a proper non--trivial submodule 
of $U$ and
consider the short exact sequence,
\begin{equation*}
0\to U_1\to U\to U_2 \to 0
\end{equation*}
Since  $\text{Ext}^1_{L(\lie g)}(U_j,V(\bpi_0)) = 0$ for $j=1,2$ by the 
induction hypothesis,
the result follows by using
the exact sequence $\text{Ext}^1_{L(\lie g)}(U_2,V(\bpi_0))\to 
\text{Ext}^1_{L(\lie g)}(U,V(\bpi_0))\to
\text{Ext}^1_{L(\lie g)}(U_1,V(\bpi_0))$. Part (ii) is now immediate by 
using a similar induction on the
length of $V_2$.
\end{pf}

The proof of proposition \ref{indecomp} is now completed as follows.
Let $V$ be an indecomposable $L(\lie g)$-module. We prove that
there exists $\chi\in\Xi$ such that $V\in\cal C_\chi$  by an  induction on 
the length of $V$.
If $V$ is irreducible, follows from the definition of
spectral characters that $V\in\cal C_{\chi_\bpi}$ for some $\bpi\in\cal P$.
If $V$ is reducible, let  $V(\bpi_0)$ be an irreducible subrepresentation of 
$V$
and let $U$ be the corresponding quotient. In other words, we have an 
extension
\begin{equation*}
0\to V(\bpi_0)\to V\to U \to 0
\end{equation*}
Write $U=\oplus_{j=1}^r U_j$ where each $U_j$ is indecomposable. By the 
inductive hypothesis,
there exist $\chi_j\in \Xi$ such that $U_j\in\cal C_{\chi_j}$, $1\le j\le 
r$.
Suppose that there exists $j_0$ such that $\chi_{j_0}\ne \chi_{\bpi_0}$. 
Lemma \ref{nullext} implies that 
$$\text{Ext}^1_{L(\lie g)}(U,V(\bpi_0)) \cong \oplus_{j=1}^r 
\rm{Ext}^1_{L(\lie g)}(U_j,V(\bpi_0))\cong
\oplus_{j\neq j_0} \rm{Ext}^1_{L(\lie g)}(U_j,V(\bpi_0)).$$
In other words, the exact sequence $0\to V(\bpi_0)\to V\to U \to 0$ is
equivalent to one of the form
\begin{equation*}
0\to V(\bpi_0)\to U_{j_0}\oplus V'\to U_{j_0}\bigoplus_{j\neq j_0} U_j \to 0
\end{equation*}
where
\begin{equation*}
0\to V(\bpi_0)\to V'\to \bigoplus_{j\neq j_0} U_j \to 0
\end{equation*}
is an element of $\oplus_{j\neq j_0} \rm{Ext}^1_{L(\lie g)}(U_j,V(\bpi_0))$. 
But this contradicts
the fact that $V$ is indecomposable. Hence $\chi_j=\chi_{\bpi_0}$ for all 
$1\le j\le r$ and $V\in\cal C_{\chi_{\bpi_0}}$.
\hfill\qedsymbol

\vskip 12pt

\section{Appendix} We give an alternate elementary  proof of Proposition 
\ref{shaded}.  This has
the advantage of computing the sequence $\mu_\ell$ of weights explicitly,
which  is useful in determining precisely the irreducible representations in 
each block. Further,
it also  makes  precise the algorithm for determining the
blocks  in
the quantum case studied in
\cite{EM}.
We proceed in two steps, namely,

\begin{enumerit}
\item[(i)] Let  $\mu\in P^+$.
There exists a sequence of weights $\mu_l \in P^+$, $l=0,\cdots,
m$, with $\mu_0 = \mu$, $\mu_m=\sum_{i\in I_\bullet} s_i\omega_i$,
$s_i\in\bn^+$,  satisfying: $$\text{Hom}_{\lie g}(\lie g\otimes
V(\mu_l), V(\mu_{l+1}))\ne 0,\ \ \forall\ 1\le l\le m.$$
\item[(ii)] Assume that $\mu=\sum_{i\in I_\bullet}
s_i\omega_i\in P^+$. Then, there exists a sequence of weights
$\mu_l \in P^+$, $l=0,\cdots, m$, with $\mu_0=\lambda_\gamma$,
$\mu_m = \mu$ satisfying $$\text{Hom}_{\lie g}(\lie g\otimes
V(\mu_l), V(\mu_{l+1}))\ne 0,\ \ \forall\ 0\le l\le m.$$
\end{enumerit}

We also need the following result proved in \cite{PRV}.
\begin{prop}\label{PRV}
  Suppose that $\lambda,\mu\in P^+$. Fix a non--zero element  $v_{w_0\mu}\in 
V(\mu)_{w_0\mu}$. Then
$V(\lambda)\otimes V(\mu)$ is generated as a $\lie g$--module by
the element $v_\lambda\otimes v_{w_0\mu}$ and the following
defining relations:
$$(x^+_{i})^{-w_0(\mu)(h_i)+1}\left(v_\lambda\otimes
v_{w_0\mu}\right)=0,\ \
(x^-_{i})^{\lambda(h_i)+1}\left(v_{\lambda}\otimes
v_{w_0\mu}\right)=0,\ \ \forall\ \ i\in I.$$
\end{prop}

\vskip 12pt

Assume that $\lie g$ is of type $A_n$ or $C_n$. Write
$\mu=\sum_{i=1}^nr_i\omega_i$. To prove the first step, we proceed
by induction on $k_0=\text{max}\{1\le k\le n,: r_k>0\}$, and show
that such a sequence exists  and further  that
$\mu_m=(\sum_{i=1}^nir_i)\omega_1$. Clearly induction starts when
$k_0=1$. Assume now that we know the result for all $k<k_0$. To
complete the inductive step we proceed by a further induction on
$r_{k_0}$. Defining $\mu_1=\mu+\sum_{i=1}^{k_0-1}\alpha_i$, it is
easily seen that $\mu_1\in P^+$ and, using Proposition \ref{PRV}
we have $$\text{Hom}_{\lie g}(\lie g\otimes V(\mu), V(\mu_{1}))\ne
0. $$ Since $$\mu_1=(r_1+1)\omega_1+\sum_{i=2}^{k_0-2}r_k\omega_k
+(r_{k_0-1}+1)\omega_{k_0-1}+(r_{k_0}-1)\omega_{k_0},$$ the proof
of step 1 is now immediate by the inductive hypothesis. To prove
the second step, it is enough to show that there exists a sequence
of the desired form if $\mu=k\omega_1$ and $\mu_0=r\omega_1$ are
such that $(k-r)\omega_1\in Q^+$. In the case of $C_n$ it suffices
to consider the case $k-r=2$. Noting that $2\omega_1=\theta$, we
see that by Proposition \ref{PRV} $$\text{Hom}_{\lie g}(\lie
g\otimes V(r\omega_1), V((r+2)\omega_{1}))\ne 0,$$ and the result
follows. For $A_n$, we have to consider the case when $k-r=n+1$.
Consider $\mu_{1}=\mu_0+\theta$ so that $$\text{Hom}_{\lie g}(\lie
g\otimes V(r\omega_1), V((r+1)\omega_{1}+\omega_n))\ne 0.$$ By the
first step we know that there exists a sequence $\mu_1,\cdots
,\mu_m$ with  $\mu_1=(r+1)\omega_{1}+\omega_n$ and
$\mu_m=(r+n+1)\omega_1$ with $$\text{Hom}_{\lie g}(\lie g\otimes
V(\mu_k), V(\mu_{k+1}))\ne 0$$ and the proof is now complete for
$A_n$.

\vskip 12pt

\noindent Suppose that  $\lie g$ is of type $B_n$  and
$\mu=\sum_{i=1}^n r_i\omega_i$. If $r_i=0$ for $i\ne n$ the first
step is obvious.  Otherwise, we have $r_k\ne 0$ for some $k<n$. We
prove by induction on $k_0=\text{min}\{1\le k< n: r_{k_0}\ne 0\}$
that we can find the sequence $\mu_1,\cdots ,\mu_m$ with
$\mu_m=(r_m+2\sum_{i=1}^n r_i)\omega_n$. When $k_0=n-1$, consider
$\mu_1=\mu+\alpha_n$. Then, Proposition \ref{PRV} implies that
$$\text{Hom}_{\lie g}\left(\lie g\otimes
V(r_{n-1}\omega_{n-1}+r_n\omega_n),
V((r_{n-1}-1)\omega_{n-1}+(r_n+2)\omega_n)\right)\ne 0,$$ and now
an obvious induction on $r_{n-1}$ gives the result.
Assume now that $k_0<n-1$ and that we know the result for all
$k>k_0$. We proceed by a further induction on $r_{k_0}$.
Set
$\mu_1=\mu+(\alpha_{k_0+1}+2(\alpha_{k_0+2}+\cdots+\alpha_n))$. We
now proceed as in the case of $A_n$ to complete the first step.
For the second step it suffices to prove the existence of the
sequence when $\mu=k\omega_n$ and $\mu_0=r\omega_n$ and $k-r =2$.
To do this observe that if we take $\mu_1=\mu+\theta$, then
$$\text{Hom}_{\lie g}(\lie g\otimes V(\mu), V(\mu_1))\ne 0,$$ and
the proof of the first step shows that we can connect $\mu_1$ and
$\mu$ by a sequence of the appropriate form.

\vskip 12pt

Suppose next  that $\lie g$ is of type  $D_n$ with $n$ even and
that $\mu=\sum_{i=1}^n r_i\omega_i$. If $r_i=0$, $i\ne n,n-1$
there is nothing to prove.  Otherwise, we have $r_k\ne 0$ for some
$k<n-1$. We prove by induction on $k_0=\text{min}\{1\le k< n-1:
r_{k_0}\ne 0\}$ that we can find two  sequences $\mu_1,\cdots
,\mu_m$, one where
$$\mu_m=(r_{n-1}+\sum_{j=0}^{\frac{n-4}{2}}r_{2j+1}+2\sum_{j=1}^{\frac{n-2}{2}}r_{2j})\omega_{n-1}
+(r_n+\sum_{j=0}^{\frac{n-4}{2}}r_{2j+1})\omega_n$$ and another
where,
$$\mu_m=(r_{n-1}+\sum_{j=0}^{\frac{n-4}{2}}r_{2j+1})\omega_{n-1}
+(r_n+\sum_{j=0}^{\frac{n-4}{2}}r_{2j+1}+2\sum_{j=1}^{\frac{n-2}{2}}r_{2j})\omega_n.$$
When $k_0=n-2$ take $\mu_1=\mu+\alpha_{n-1}$ (resp.
$\mu_1=\mu+\alpha_{n}$) and proceed by induction on $r_{n-2}$. To
complete the inductive step for $k_0<n-2$ we take
$\mu_1=\mu+\alpha_{k_0+1}+2(\alpha_{k_0+2}+\cdots+
\alpha_{n-2})+\alpha_{n-1}+\alpha_n$, we omit further details. For
the second step, we must prove that $k\omega_i$ and
$(k-2)\omega_i$ are connected by an appropriate sequence of
elements of $P^+$ for $i=n,n-1$. As before, we take
$\mu_1=(k-2)\omega_i+\theta$ and use the first step to get the
result.

Now consider the case of $D_n$ with $n$ odd and let
$\mu=\sum_{i=1}^n r_i\omega_i$. If $r_i=0$, $i\ne n$ there is
nothing to prove.  In the general case we  proceed in two further
steps:
\begin{enumerit}
\item[(a)]  There exists a sequence of weights $\mu_l \in P^+$, $l=0,\cdots,
m$, with $\mu_0 = \mu$, $\mu_m=\sum_{i\text{ odd}} s_i\omega_i$,
$s_i\in\bn^+$,  satisfying: $$\text{Hom}_{\lie g}(\lie g\otimes
V(\mu_l), V(\mu_{l+1}))\ne 0,\ \ \forall\ 1\le l\le m.$$

\item[(b)]  Assume that $\mu$ is supported only on the odd nodes. Then, 
there exists a sequence of weights $\mu_l\in P^+$, $l=0,\cdots,
m$, with $\mu_0 = \mu$, $\mu_m=\sum_{i\in I_\bullet} s_i\omega_i$,
$s_i\in\bn^+$,  satisfying: $$\text{Hom}_{\lie g}(\lie g\otimes
V(\mu_l), V(\mu_{l+1}))\ne 0,\ \ \forall\ 1\le l\le m.$$
\end{enumerit}
To prove step (a) we assume that $r_k>0$ for some $k$ even and
proceed by induction on $k_0 = {\rm{min}}\{k \text{ even}:
r_k>0\}$. First assume that  $k_0 = n-1$ and proceed by a further
induction on $r_{n-1}$ as usual. Setting
$$\mu_1=\mu+(\alpha_1+\cdots +\alpha_{n-2}+\alpha_n) =
\sum_{j=1}^{\frac{n-3}{2}}r_{2j+1}\omega_{2j+1}
+ (r_1+1)\omega_1+(r_n+1)\omega_n+(r_{n-1}-1)\omega_{n-1},$$
and using the induction on $r_{n-1}-1$ completes this case. Next,
  suppose that $k_0 = n-3$ and
take  $$\mu_1 = \mu+(\alpha_{n-2}+\alpha_{n-1}+\alpha_{n}) =
  \sum_{j=0}^{\frac{n-3}{2}}r_{2j+1}\omega_{2j+1} + 
(r_{n-1}+1)\omega_{n-1}+(r_n+1)\omega_n+(r_{n-3}-1)\omega_{n-3}$$
   and the result follows by
   induction on $r_{n-3}$.
   Now assume that $k_0<n-3$ and that we know the result for all $k>k_0$.
   Taking  $\mu_1 = 
\mu+(\alpha_{k_0+1}+2(\sum_{i=k_0+2}^{n-2}\alpha_i)+\alpha_{n-1}+\alpha_{n-2})$. 
Then
$$\mu_1 = \sum_{j=0}^{\frac{n-1}{2}}r_{2j+1}\omega_{2j+1} +
\sum_{j=\frac{k_0+4}{2}}^{\frac{n-1}{2}}r_{2j}
\omega_{2j}+(r_{k_0+2}+1)\omega_{k_0+2}+(r_{k_0}-1)\omega_{k_0}$$
completes the inductive step. Observe that when $k_0 = 2$ we have
$$\mu_m = \sum_{j=1}^{\frac{n-3}{2}}r_{2j+1}\omega_{2j+1} +
(r_1+\sum_{j=1}^{\frac{n-1}{2}} r_{2j}) \omega_1 +
(r_n+r_{n-1}+2\sum_{j=1}^{\frac{n-3}{2}} r_{2j})\omega_n$$

Now we prove step (b), i.e., $r_j=0$ for all $1\le j\le n$ with $j$
even.  We proceed by induction on $k_0 = {\rm{min}}\{k ; r_k>0\}$
and on $r_{k_0}$. If $k_0=n$ there is nothing to prove. If $k_0 =
n-2$, then taking $\mu_1 = \mu+\alpha_n =
(r_n+2)\omega_n+(r_{n-2}-1)\omega_{n-2}$  completes the induction.
Now assume that $k_0<n-2$. Taking $\mu_1 =
\mu+(\alpha_{k_0+1}+2(\sum_{i=k_0+2}^{n-2}\alpha_i)+\alpha_{n-1}+\alpha_{n-2})$
we  see that $$\mu_1 =
\sum_{j=\frac{k_0+3}{2}}^{\frac{n-1}{2}}r_{2j+1}\omega_{2j+1} +
(r_{k_0+2}+1)\omega_{k_0+2}+(r_{k_0}-1)\omega_{k_0}$$ This
completes the proof of the first step, notice that following this
procedure gives $$\mu_m
=\left(r_n+3r_{n-1}+2\sum\limits_{j=0}^{\frac{n-3}{2}}r_{2j+1}+4\sum\limits_{j=1}^{\frac{n-3}{2}}r_{2j}\right)\omega_n$$

The second step is completed by the usual method and we omit all
details.

\vskip 12pt

\noindent  {\bf{$\lie g=E_6$}} Consider the following sequence of
weights:
\begin{eqnarray*}
&\lambda_1=&(r_1+r_6)\omega_1+r_2\omega_2+r_3\omega_3+r_4\omega_4+(r_5+r_6)\omega_5,\\
&\lambda_2=&(r_1+r_3+r_6)\omega_1+(r_2+r_3)\omega_2+r_4\omega_4+(r_5+r_6)\omega_5,\\
&\lambda_3=&
(r_1+r_3+r_6)\omega_1+(r_2+r_3)\omega_2+(2r_4+r_5+r_6)\omega_5,\\
&\lambda_4=&
(r_1+r_3+r_6)\omega_1+(r_2+r_3+2r_4+r_5+r_6)\omega_2,\\
&\lambda_5=&(r_1+2r_2+3r_3+4r_4+2r_5+3r_6)\omega_1.\end{eqnarray*}
Setting $\mu=\lambda_0$, it suffices to show that $\lambda_k$ and
$\lambda_{k+1}$ are connected by a sequence of weights as in (i)
above. But this is clear from Proposition \ref{PRV}, by noting that
\begin{eqnarray*}&\lambda_1-\lambda_0&=r_6(\alpha_1+\alpha_2+\alpha_3+\alpha_4+\alpha_5),\\
&\lambda_2-\lambda_1&=r_3(\alpha_1+\alpha_2),\\
&\lambda_3-\lambda_2&=r_4\alpha_5,\\
&\lambda_4-\lambda_5&=(2r_4+r_5+r_6)(\alpha_1+2\alpha_2+2\alpha_3+\alpha_4+\alpha_6),\\
&\lambda_5-\lambda_6&=(r_2+r_3+2r_4+r_5+r_6)\alpha_1.\end{eqnarray*}

To prove the second step we can assume that $\mu_0=r\omega_1, \mu
= k\omega_1$ and $k-r = 3$. Take
$$\mu_1=\mu_0+\theta=\mu+\omega_6.$$ Then by  Proposition \ref{PRV},
we have $\text{Hom}_{\lie g}(\lie g\otimes V(\mu), V(\mu_{1}))\ne
0$. On the other hand, we see from step (i) that there exists an
appropriate sequence connecting $\mu_1$ and $(r+3)\omega_1$. The
result is proved for $E_6$.

\vskip 12pt

{\bf{$\lie g=E_7$}}  Consider the following sequence of weights:
\begin{eqnarray*}
&\lambda_1=&(r_1+r_7)\omega_1+r_2\omega_2+r_3\omega_3+r_4\omega_4+r_5\omega_5+(r_6+r_7)\omega_6,\\
&\lambda_2=&(r_1+r_4+r_7)\omega_1+r_2\omega_2+(r_3+r_4)\omega_3+r_5\omega_5+(r_6+r_7)\omega_6,\\
&\lambda_3=&
(r_1+r_4+r_7)\omega_1+r_2\omega_2+(r_3+r_4)\omega_3+(r_6+r_7+2r_5)\omega_6,\\
&\lambda_4=&(r_1+r_4+r_7)\omega_1+(r_2+r_6+r_7+2r_5)\omega_2+(r_3+r_4)\omega_3
\\
&\lambda_5=&(r_1+r_3+2r_4+r_7)\omega_1+(r_2+r_3+r_4+2r_5+r_6+r_7)\omega_2,\\
&\lambda_6=&(r_1+2r_2+3r_3+4r_4+4r_5+2r_6+3r_7)\omega_1.\end{eqnarray*}
Setting $\mu=\lambda_0$ we see again that $\lambda_k$ and
$\lambda_{k+1}$ are connected by an appropriate sequence. For the
second step, we can assume that $\mu_0 = r\omega_1, \mu =
k\omega_1$ with $k-r = 2$. Taking $\mu_1=\mu_0+\theta=\mu+\omega_6$,
we find from step (i) that $\mu_1$ and $(k+2)\omega_1$ are
connected and we are done.

\vskip 12pt

{\bf{$\lie g=E_8$}}  Consider the following sequence of weights:
\begin{eqnarray*}
&\lambda_1=&(r_1+r_8)\omega_1+r_2\omega_2+r_3\omega_3+r_4\omega_4+r_5\omega_5+r_6\omega_6+(r_7+r_8)\omega_7,\\
&\lambda_2=&(r_1+r_5+r_8)\omega_1+r_2\omega_2+r_3\omega_3+(r_4+r_5)\omega_4+r_6\omega_6+(r_7+r_8)\omega_7,\\
&\lambda_3=&
(r_1+r_5+r_8)\omega_1+r_2\omega_2+r_3\omega_3+(r_4+r_5)\omega_4+(r_7+r_8+2r_6)\omega_7,\\
&\lambda_4=&(r_1+r_4+r_7)\omega_1+(r_2+2r_6+r_7+r_8)\omega_2+r_3\omega_3+(r_4+r_5)\omega_4
\\
&\lambda_5=&(r_1+2r_4+r_5+r_7)\omega_1+(r_2+2r_6+r_7+r_8)\omega_2+(r_3+r_4+r_5)\omega_3,
\\
&\lambda_6=&(r_1+r_3+3r_4+2r_5+r_7)\omega_1+(r_2+r_3+r_4+r_5+2r_6+r_7+r_8)\omega_2,\\
&\lambda_7=&(r_1+2r_2+3r_3+5r_4+4r_5+4r_6+3r_7+2r_8)\omega_1
.\end{eqnarray*} Setting $\mu=\lambda_0$ we see again that
$\lambda_k$ and $\lambda_{k+1}$ are connected by an appropriate
sequence. For the second step, we can assume that $\mu_0 =
r\omega_1, \mu = k\omega_1$ with $k-r = 1$. Taking
$\mu_1=\mu_0+\theta=\mu_0+\omega_1 = \mu$ and we are done.

\vskip 12pt

{\bf{$\lie g=F_4$}} Consider the following sequence of weights:
\begin{eqnarray*}
&\lambda_1&=(r_1+2r_2)\omega_1+r_3\omega_3+r_4\omega_4,\\
&\lambda_2&=(r_1+2r_2)\omega_1+(r_4+2r_3)\omega_4,\\
&\lambda_3&=(r_1+2r_2+4r_3+2r_4)\omega_1.
\end{eqnarray*}
Setting $\mu=\lambda_0$ we see again that $\lambda_k$ and
$\lambda_{k+1}$ are connected by an appropriate sequence.
For the second step we can assume that $\mu = r\omega_1$, with $r\neq 0$. 
Then we define,
Then set \begin{eqnarray*} &\mu_1 =&
\mu+(\alpha_1+3\alpha_2+2\alpha_3+\alpha_4),\\ & \mu_2 =& \mu_1
+\alpha_1, \\ & \mu_3 =& \mu_2 -(2\alpha_1+2\alpha_2+\alpha_3),\\ & \mu_4 =& 
\mu_3 -
\theta\end{eqnarray*} and the result is proved by induction on $r$, noting 
that $\mu_4 = (r-1)\omega_1$.

{\bf{$\lie g=G_2$}} Here we define
$\lambda_1=\mu+r_2(3\alpha_1+\alpha_2)$ to see that $\mu$ and
$(r_1+3r_2)\omega_1$ are connected as in Step (i). To prove step
(ii), we use the fact that
$r\omega_1+(2\alpha_1+\alpha_2)=(r+1)\omega_2$ to get the result.

\bibliographystyle{amsplain}

\begin{thebibliography}{10}

\bibitem{AK}
T.~Akasaka and M.~Kashiwara, {\em Finite-dimensional
representations of
quantum affine algebras}, Publ. Res. Inst. Math. Sci. \textbf{33} (1997),
no.~5, 839--867.

\bibitem{C} V.~Chari, {\em Integrable representations of affine
{L}ie algebras}, Invent.
Math. {\bf 85} (1986), no.~2, 317--335.

\bibitem{C1} \bysame,
{\em Braid group actions and tensor products.} Internat. Math.
Res. Notices (2002), no. 7, 357--382.

\bibitem{CPnew}
V.~Chari and A.~Pressley, {\em New unitary representations of loop
groups},
Math. Ann. {\bf 275} (1986), no.~1, 87--104.
\bibitem{CPqa}
\bysame, {\em Quantum affine algebras}, Comm. Math. Phys. {\bf
142} (1991),
no.~2, 261--283.

\bibitem{CPweyl}
\bysame, {\em Weyl modules for classical and quantum affine
algebras.} Represent. Theory {\bf 5} (2001), 191--223

\bibitem{D}
J. Dixmier. {\em Algebres Envelopantes}, Cahiers Scientifiques {\bf 37}, (1974).

\bibitem{EM}
P.~Etingof and A.~Moura. {\em Elliptic Central Characters and
Blocks of Finite Dimensional Representations of Quantum Affine
Algebras}, {\em Represent. Theo.} {\bf 7}, (2003), 346--373.

\bibitem{FL1}B.~Feigin and S.~Loktev,  On generalized Kostka polynomials and 
quantum Verlinde
rule, To appear in D.Fuchs 60-th anniversary volume, Preprint,
math.QA/9812093


\bibitem{FL2}\bysame,  Multi-dimensional Weyl Modules and Symmetric 
Functions, Preprint, math.QA/0212001.

\bibitem{FM}
E. Frenkel\ and\ E. Mukhin, {\em Combinatorics of $q$-characters
of finite-dimensional representations of quantum affine algebras},
Comm. Math. Phys. {\bfseries 216} (2001), no.~1, 23--57

\bibitem{FR}
E.~Frenkel and N.~Reshetikhin, {\em The $q$-characters of
representations of
quantum affine algebras and deformations of $\mathcal{W}$-algebras}, Recent
developments in quantum affine algebras and related topics (Raleigh, NC,
1998), Contemp. Math., 248, Amer. Math. Soc., Providence, RI, 1999,
pp.~163--205.

\bibitem{H} J.~E.~Humphreys,
{\em Introduction to Lie Algebras and Representation
Theory}, Graduate Texts in Mathematics, Springer--Verlag, 1970.


\bibitem{Ka2}
M.~Kashiwara, {\em On level-zero representation of quantized
affine algebras.}, Duke Math. J. {\bf 112} (2002), no.~1, 117--195

\bibitem{K}
B. Kostant, {\em Lie Group representations on polyomial rings.}, Am.J.Math. {\bf 85} (1963),  327--404.



\bibitem{N1}
H. Nakajima, {\em $t$-analogue of the $q$-characters of finite
dimensional representations of quantum affine algebras}, Physics
and combinatorics, 2000 (Nagoya), 196--219, World Sci. Publishing,
River Edge, NJ, 2001

\bibitem{N2}
\bysame, {\em Extremal weight modules of quantum affine algebras},
Preprint math.QA/0204183.

\bibitem{PRV} K.R.Parthasarathy, R. Ranga Rao, and
V.S.Varadarajan, {\em Representations of complex semi--simple Lie
groups and Lie Algebras}, Ann.Math. {\bf 85}, (1967), 38--429.

\end{thebibliography}

\end{document}